# Notes on Symmetric Bases

**History**

> 0.01  17-02-93   First version completed.
> 0.02  22-09-14   Edited for ArxIv publication. New Abstract and Reference sections added, and all references in the body of the paper updated to reflect the new reference numbers. No other changes made.

**Abstract**

$A_k = \{1, a_2, \ldots a_k\}$ is an h-basis for n if every positive integer not exceeding n can be expressed as the sum of no more than h values $a_i$; we write $n = n_h(A_k)$. An extremal h-basis $A_k$ is one for which n is as large as possible. Computing such extremal bases has become known as the Postage Stamp Problem.

A basis $A_k$ is *symmetric* if $A_k = \{1, a_2, \ldots a_k\}$ where $a_i + a_{k-i} = a_k$ for $1 <= i <= k-1$. Examination of a number of symmetric bases suggests the following conjecture: if the range $0 \ldots a_k$ is covered using at most h stamps, then the range $0 \ldots ha_k$ is also covered using at most h stamps. This paper shows that this is not strictly true, but demonstrates that there is a value $h_1$ such that the conjecture is true for all $h >= h_1$.

Some of the content of the paper is derived directly from Selmer's monographs (see, for example, Selmer, E.S., [6] page 8.12), although the proof of the special case of Meure's theorem is my own.

## 1  Background

*1.1  Introduction*

A basis $A_k$ is *symmetric* if:

$$A_k = \{1, a_2, \ldots a_k\} \quad \text{where} \quad a_i + a_{k-i} = a_k \quad \text{for } 1<=i<=k-1 \quad - (1)$$

Examination of a number of symmetric bases suggests the following conjecture:

> If the range $0 \ldots a_k$ is covered using at most h stamps, then the range $0 \ldots ha_k$ is also covered using at most h stamps.

In fact, this is not quite true; however, it can be proved that there is a value $h_1$ such that the conjecture is true for all $h>=h_1$.

The following notes are derived directly from Selmer's monographs (Selmer, E.S., [6], [7], [8]) although the proof of the special case of Meure's theorem is my own.

*1.2  Admissibility and the definition of $h_0$*

Given any basis $A_k$, we can look at the cover $n(h, A_k)$ for increasing values of h starting from 1; for example, take $A_k = \{1, 3, 6, 10\}$:

$n(1, A_k) = 1$
$n(2, A_k) = 4$
$n(3, A_k) = 23$

We say that a basis is *h-admissible* only if its cover exceeds $a_k$, and we define $h_0$ to be the smallest value of h for which this is true; in the example above, $h_0 = 3$.

We can now rephrase the *conjecture* as follows:

For any symmetric basis $A_k$, $n(h, A_k) = ha_k$ for all $h >= h_0$

and the *theorem* becomes:

For any symmetric basis $A_k$ there exists a value $h_1 >= h_0$ such that $n(h, A_k) = ha_k$ for all $h >= h_1$.

## 2  Counter-examples to the conjecture

It has been proved that there are no counter-examples to the conjecture for k<=6, and no counter-example has been found for k=7 or k=8 (but neither has it been proved that none exist); see (Selmer, E.S., [6] page 8.12).

The "simplest" counter-example known is for k=9 as follows:

$A_9 = \{1, 3, 5, 8, 20, 23, 25, 27, 28\}$

$n(2, A_9) = 6$     ( $< a_9$ )
$n(3, A_9) = 41$     ( $> a_9, < 3a_9$ )
$n(4, A_9) = 112$     ( $= 4a_9$ )

Here $h_0 = 3$ and $h_1 = 4$.

This counter-example is, in fact, one of a family of such bases derived from the following parametric basis $A_5$:

$A_5(p) = \{ 1, p, p+2, 2p+2, (3p^2+3p+4)/2 \}$    for p odd, p>=3

The symmetric basis $A_9(p)$ is created by simply extending the basis $A_5(p)$ in accord with the rules of symmetry; eg for p = 3 we have:

$A_5(3) = \{1, 3, 5, 8, 20\}$ with differences $\{1, 2, 2, 3, 12\}$

which, when extended, gives us the basis $A_9$ above with differences $\{1, 2, 2, 3, 12, 3, 2, 2, 1\}$.

The basis $A_5(p)$ can also be extended in the obvious way to produce a symmetric basis $A_{10}(p)$, and for p odd, p>=5, this parametric basis, too, is a counter-example to the conjecture; eg:

$A_5(5) = \{1, 5, 7, 12, 47\}$ with differences $\{1, 4, 2, 5, 35\}$

can be extended to:

$A_{10}(5) = \{1, 5, 7, 12, 47, 82, 87, 89, 93, 94\}$ with differences $\{1, 4, 2, 5, 35, 35, 5, 2, 4, 1\}$

and we find:

$n(4, A_{10}) = 34$    $(< a_{10})$
$n(5, A_{10}) = 132$    $(> a_{10}, < 5a_{10})$
$n(6, A_{10}) = 564$    $(= 6a_{10})$

thus showing that $h_0 = 5, h_1 = 6$.

In fact, for bases $A_9(p)$ and $A_{10}(p)$ derived as shown above, it can be proved that $h_0=p, h_1=h_0+1$.

## 3  Proof of the theorem

*3.1 Lemma*

We first prove the following lemma:

If $0 <= x < a_k$, then $h_0 a_k - x$ has an $h_0$-generation.

Proof:

$x$ is less than $a_k$, and so, by definition of $h_0$, has an $h_0$-generation, say:

$$x = c_{k-1}a_{k-1} + c_{k-2}a_{k-2} + ... + c_2 a_2 + c_1 a_1 \qquad \text{with} \quad c_{k-1}+c_{k-2}+ ... +c_1 <= h_0 \qquad - (2)$$

Because the basis $A_k$ is symmetric, we can use (1) to re-write this as:

$$x = c_{k-1}(a_k - a_1) + c_{k-2}(a_k - a_2) + ... c_2(a_k - a_{k-2}) + c_1(a_k - a_{k-1})$$

and hence

$$(c_{k-1}+c_{k-2}+ ... +c_1)a_k - x = c_{k-1}a_1 + c_{k-2}a_2 + ...+ c_2 a_{k-2} + c_1 a_{k-1}$$

Let $c_k = h_0 - (c_{k-1}+c_{k-2}+ ... +c_1)$; clearly $c_k >= 0$, and we can add $c_k a_k$ to both sides:

$$h_0 a_k - x = c_k a_k + c_1 a_{k-1} + c_2 a_{k-2} + ... + c_{k-2} a_2 + c_{k-1} a_1 \quad \text{with} \quad c_k + c_{k-1} + c_{k-2} + ... + c_1 = h_0$$

Thus our lemma is proven: the line above is an $h_0$-generation of $h_0 a_k - x$.

*3.2 The theorem*

Choose $h = 2h_0 - 2$.

Then the lemma shows that all values:

$$0 \leq x \leq a_k \quad \text{and} \quad (h_0-1)a_k \leq x \leq h_0 a_k \quad - (3)$$

have $h_0$-generations, and so all values:

$$k a_k \leq x \leq (k+1)a_k \quad \text{and} \quad (k+h_0-1)a_k \leq x \leq (k+h_0)a_k \quad - (4)$$

have h-generations, provided that $1 \leq k \leq h_0-2$ [since up to $h-h_0=h_0-2$ stamps $a_k$ can be added to each of the generations needed for (3) above].

But the ranges (3) and (4) for $1 \leq k \leq h_0-2$ are contiguous and together form the range:

$$0 \leq x \leq (h_0-2+h_0)a_k = h a_k$$

and so we have $n(h, A_k) = h a_k$.

This shows that there must exist some minimal value $h_1$ such that:

$$n(h_1, A_k) = h_1 a_k \quad \text{and} \quad h_0 \leq h_1 \leq 2h_0-2$$

and so the theorem is proved.

## 4 Remark

In Selmer's monograph, the theorem above is deduced from a more general formula called *Meure's formula*, which makes a connection between the Frobenius "coin" problem and the postage stamp problem. The statement and proof of Meure's formula is more complex than the simple proof given above, but furnishes a more powerful result as follows:

> If $A_k = \{1, a_2, \ldots a_k\}$ is a basis such that $a_{k-1} = a_k-1$, then there exists an $h_1 \geq h_0$ such that $n(h, A_k) = h a_k$ for all $h \geq h_1$.

See (Selmer, E.S., [6] page 7.7).